\documentclass[dvips,epsfig,reqno,graphics]{amsart}
\usepackage{amssymb}
\usepackage{amsfonts}
\usepackage{amsmath}
\usepackage{amsxtra}

\title[Translation Weingarten Surfaces]
  {Polynomial Translation Weingarten Surfaces in $3$-dimensional Euclidean space}
\author{Marian Ioan Munteanu}
\author{Ana Irina Nistor}

\newtheorem{contor}{1.}

\newtheorem{theorem}[contor]{Theorem}

\def\proof{{\sc Proof.\ }}

\def\R{\mathbb{R}}
\def\E{\mathbb{E}}
\def\a{\alpha}

\def\b{\beta}

\date{ \today }

\begin{document}

\maketitle

\input epsfx.tex

\begin{abstract}
In this paper we will classify those translation surfaces in $\E^3$
involving polynomials which are Weingarten surfaces.

\vspace{1mm}


\vspace{1mm}


\end{abstract}

\section{ Preliminaries}

The surfaces of constant mean curvature, $H$-surfaces and those of constant
Gaussian curvature, $K$-surfaces in the Euclidean 3-dimensional space, $\E^3$,
have been studied extensively.
One interesting class of surfaces in $\E^3$ is that of {\it translation surfaces},
which can be parametrized, locally, as $r(u,v)=(u,v,f(u)+g(v))$, where $f$ and $g$ are smooth
functions. This type of surfaces are important either because they are interesting
themselves or because they furnish counterexamples for some problems
(e.g. it is a known fact that a minimal surface has vanishing second Gaussian curvature
but not conversely -- see for details \cite{kn:BK92}).
We call {\em polynomial translation surfaces} (in short, {\em PT surfaces})
those translation
surfaces for which $f$ and $g$ are polynomials.

Scherk's surface, obtained by H. Scherk in 1834, is the only non flat minimal surface,
that can be represented as a translation surface. More precisely we have
{\bf Theorem A.}
{\em Let $S$ be a translation minimal surface in $3$-dimensional Euclidean space.
Then $S$ is an open part of $\E^3$ or it is congruent to the following surface
$$
   z=\frac 1a\ \log \left| \frac{\cos(a x)}{\cos(a y)} \right|,\quad a\neq 0.
$$
}

Other interesting results concerning translation surfaces having
either constant mean curvature or constant Gaussian curvature are
the following:

{\bf Theorem B.} {\em Let $S$ be a translation surface with constant Gauss curvature $K$ in
$3$-dimensional Euclidean space. Then $S$ is congruent to a cylinder, so it is flat}
($K=0$) (Theorem 1 in \cite{kn:Liu99}).

{\bf Theorem C.} {\em Let $S$ be a translation surface with constant mean curvature $H\neq0$
in $3$-dimensional Euclidean space $\E^3$. Then $S$ is congruent to the following
surface}
$$
     z=\frac{\sqrt{1+a^2}}{2H}\ \sqrt{1-4H^2x^2}+a y,\quad a\in\R
$$
(Theorem 2, statement (1) in \cite{kn:Liu99}).
Cf. also References in \cite{kn:Liu99}.

\medskip

A surface $S$ is called a {\em Weingarten surface} if there is some (smooth)
relation $W(\kappa_1,\kappa_2)=0$ between its two principal curvatures $\kappa_1$
and $\kappa_2$, or equivalently, if there is a (smooth) relation $U(K,H)=0$ between its
mean curvature $H$ and its Gaussian curvature $K$.

\medskip

In this paper we study those PT surfaces that are Weingarten surfaces.
We abbreviate by a {\em WPT surfaces}.

We give the following classification theorem:

\begin{theorem}
Let $S$ be a WPT surface in $3$-dimensional Euclidean space.
Then
\begin{itemize}
\item [(i)] $S$ is a cylinder, case in which $K=0$;
\item [(ii)] $S$ is a paraboloid of revolution, case in which the mean curvature $H$
and the Gaussian curvature $K$ are positive everywhere and related by the relation
\begin{equation}
\label{eq:cond_parab}
   8aH^2=\sqrt{K}(2a+\sqrt{K})^2
\end{equation}
where $a$ is a positive constant.
\end{itemize}
\end{theorem}

\section{Weingarten translation surfaces}

Let $r:S\longrightarrow\R^3$ be an isometrical immersion of a translation surface
of type
\begin{equation}
\label{eq:param}
r(u,v)=(u,v,f(u)+g(v))
\end{equation}
where $f$ and $g$ are smooth functions.
The first fundamental form ${\bf{I}}$
and the second fundamental form ${\bf II}$
have particular forms, namely
$$
{\bf{I}}=\left(1+f'(u)^2\right)du^2+2 f'(u)g'(v)du\ dv+\left(1+g'(v)^2\right)dv^2
$$
$$
{\bf II}=\frac1{\sqrt\Delta}\ \big(f''(u)\ du^2+g''(v)\ dv^2\big)
$$
where $\Delta=1+f'(u)^2+g'(v)^2$.
Let us denote $f'$ by $\a$ and $g'$ by $\b$.
Hence, the mean curvature $H$ and the Gaussian curvature $K$ can be written as
\begin{equation}
\label{eq:H}
H=\frac{\left(1+\beta(v)^2\right)\alpha'(u)+\left(1+\alpha(u)^2\right)\beta'(v)}
       {2\big[1+\alpha(u)^2+\beta(v)^2\big]^{\frac32}}
\end{equation}
\begin{equation}
\label{eq:K}
K=\frac{\alpha'(u)\beta'(v)}{\big[1+\alpha(u)^2+\beta(v)^2\big]^2}\ .
\end{equation}

The existence of a Weingarten relation $U(H,K)=0$ means that curvatures
$H$ and $K$ are functionally related, and since $H$ and $K$ are
differentiable functions depending on $u$ and $v$,
this implies the {\em Jacobian condition}
$\frac{\partial(H,K)}{\partial(u,v)}=0$.
More precisely the following condition
\begin{equation}
\frac{\partial H}{\partial u}\frac{\partial K}{\partial v}-
\frac{\partial H}{\partial v}\frac{\partial K}{\partial u}=0
\end{equation}
needs to be satisfied.
Conversely, if the above condition is satisfied, then the
curvatures must be functionally related and thus, by definition, the surface
must be Weingarten. The Jacobian condition characterizes W-surfaces and
it is used to identify them when an explicit Weingarten relation
cannot be found.

In our case, the Jacobian condition yields the following relation
\begin{equation}
\label{eq:Jac_cond}
\begin{array}{c}
 8\alpha(u)\beta(v)\alpha'(u)^3\beta'(v)^2-
 8\alpha(u)\beta(v)\alpha'(u)^2\beta'(v)^3-\\
- 3\beta(v)\alpha'(u)\beta'(v)^2\alpha''(u)+
 3\alpha(u)\alpha'(u)^2\beta'(v)\beta''(v)-\\
-2\alpha(u)^2\beta(v)\alpha'(u)\beta'(v)^2\alpha''(u)+
 2\alpha(u)\beta(v)^2\alpha'(u)^2\beta'(v)\beta''(v)-\\
-3\beta(v)^3\alpha'(u)\beta'(v)^2\alpha''(u)+
 3\alpha(u)^3\alpha'(u)^2\beta'(v)\beta''(v)-\\
-3\alpha(u)\alpha'(u)^3\beta''(v)+
 3\beta(v)\beta'(v)^3\alpha''(u)+\\
+3\alpha(u)^2\beta(v)\beta'(v)^3\alpha''(u)-
 3\alpha(u)\beta(v)^2\alpha'(u)^3\beta''(v)+\\
+\alpha'(u)\alpha''(u)\beta''(v)-
 \beta'(v)\alpha''(u)\beta''(v)+\\
+ \alpha(u)^2\alpha'(u)\alpha''(u)\beta''(v)-
 \beta(v)^2\beta'(v)\alpha''(u)\beta''(v)+\\
+2\beta(v)^2\alpha'(u)\alpha''(u)\beta''(v)-
 2\alpha(u)^2\beta'(v)\alpha''(u)\beta''(v)+\\
+\alpha(u)^2\beta(v)^2\alpha'(u)\alpha''(u)\beta''(v)-
 \alpha(u)^2\beta(v)^2\beta'(v)\alpha''(u)\beta''(v)-\\
-\alpha(u)^4\beta'(v)\alpha''(u)\beta''(v)+
 \beta(v)^4\alpha'(u)\alpha''(u)\beta''(v)=0.
\end{array}
\end{equation}

At this point we will consider $\alpha$ and $\beta$ to be polynomials
of degree $m$ and $n$ respectively. More precisely we shall consider
$$
\alpha=a_mu^m+a_{m-1}u^{m-1}+\ldots\quad {\rm and}\quad
\beta=b_nv^n+b_{n-1}v^{n-1}+\ldots
$$
where $a_m$ and $b_n$ are different from $0$.
Replacing $\alpha$ and $\beta$ in \eqref{eq:Jac_cond} we obtain
a polynomial expression in $u$ and $v$ vanishing identically.
This means that all the coefficients are 0.

\bigskip

Let us distinguish several cases:

\medskip

{\bf Case 1:} $m,n\geq 2$, i.e. $\alpha''\neq 0$ and $\beta''\neq 0$.

\medskip

{\bf a.} Suppose $m>n (\geq2)$

The dominant term corresponds to $u^{5m-2}v^{2n-3}$ and it comes from
$$
  3\alpha(u)^3\alpha'(u)^2\beta'(v)\beta''(v)-\alpha(u)^4\beta'(v)\alpha''(u)\beta''(v)
$$
having the coefficient
$$
 a_m^5b_n^2\left(3m^2n^2(n-1)-mn^2(m-1)(n-1)\right).
$$
This cannot vanish since $a_m,b_n\neq 0$ and $m>n\geq2$.

\smallskip

The subcase $n>m\geq2$ can be treated in similar way.

\medskip

{\bf b.} Suppose $m=n\geq2$

In the same manner, this case cannot occur.

\bigskip

{\bf Case 2:} $m>n=1$

In this case $\beta$ can be expressed as $\beta(v)=a v + b$, with $a$ and $b$
real constants, $a\neq 0$. We rewrite the Jacobian condition in the following way
\begin{equation}
\label{eq:Jac_cond_n1}
\begin{array}{c}
 8a^2\alpha(u)\beta(v)\alpha'(u)^3-
 8a^3\alpha(u)\beta(v)\alpha'(u)^2-
 3a^2\beta(v)\alpha'(u)\alpha''(u)-\\
-2a\alpha(u)^2\beta(v)\alpha'(u)\alpha''(u)-
 3a^2\beta(v)^3\alpha'(u)\alpha''(u)+\\
+ 3a^3\beta(v)\alpha''(u)+
 3a^3\alpha(u)^2\beta(v)\alpha''(u)=0.
\end{array}
\end{equation}

Let us analyze the terms in $u$ having maximum degree, namely $u^{4m-3}v$.
This comes from $-2a\alpha(u)^2\beta(v)\alpha'(u)\alpha''(u)$ and yields
the relation $2a^2a_m^4m^2(m-1)=0$ which cannot hold in this case.

\smallskip

The case $n>m=1$ can be treated in similar way.

\medskip

{\bf Case 3:} $m=n=1$

In this case $\alpha$ and $\beta$ can be expressed as $\alpha(u)=Au+B$ and
$\beta(v)=av+b$, with $A$, $B$, $a$ and $b$ real constants, $A,a\neq0$.
The Jacobian condition becomes
\begin{equation}
\label{eq:Jac_cond_mn1}
8\alpha(u)\beta(v)\alpha'(u)^3\beta'(v)^2-
8\alpha(u)\beta(v)\alpha'(u)^2\beta'(v)^3=0.
\end{equation}
Using the same technique as above one gets $A=a$.
So, the parametrization of the surface can be written
(after a possible translation in $\E^3$) in the form
\begin{equation}
\label{eq:surf_mn1}
r(u,v)=\big( u,v,a(u-u_0)^2+a(v-v_0)^2 \big)
\end{equation}
where $u_0,v_0\in\R$.
This is a paraboloid of revolution, and its curvatures $H$ and $K$
are both everywhere positive and they are related by the relation
\eqref{eq:cond_parab}.

\begin{figure}[htb]
\begin{center}
\epsfxsize=50mm
\centerline{\leavevmode
\epsffile{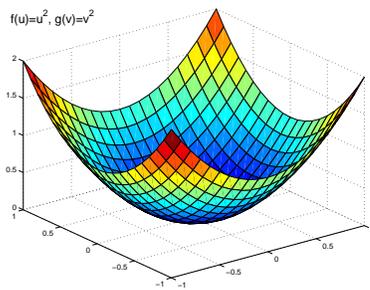}
}
\end{center}
\caption{The paraboloid of revolution}
\end{figure}

\medskip

{\bf Case 4:} $m=0$ (or $n=0$)

In this case $\alpha$ (or $\beta$) is constant and the Jacobian condition
is automatically satisfied. So, the parametrization of the surface
can be written in the form
\begin{eqnarray}
\label{eq:surf_m0}
r(u,v)=\big( u,v, au+g(v) \big)\\
\label{eq:surf_n0}
r(u,v)=\big( u,v, f(u)+a v \big)
\end{eqnarray}
where $f$ and $g$ are arbitrary polynomials, $a\in\R$ (it can also vanishes).
These two surfaces are cylinders and they are obviously flat.

\begin{figure}[htb]
\begin{center}
\epsfxsize=50mm
\centerline{\leavevmode
\epsffile{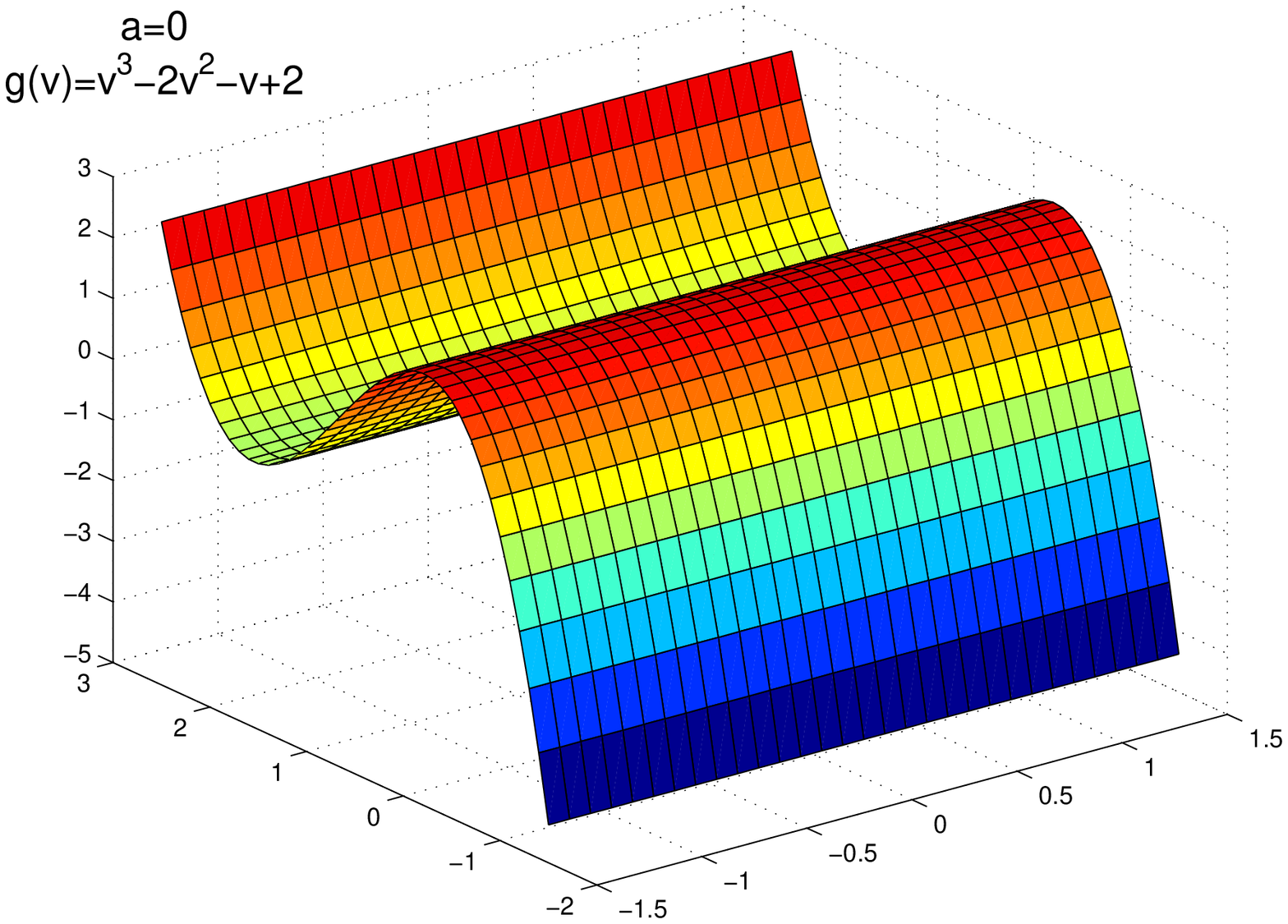}
\qquad
\epsfxsize=50mm
\epsffile{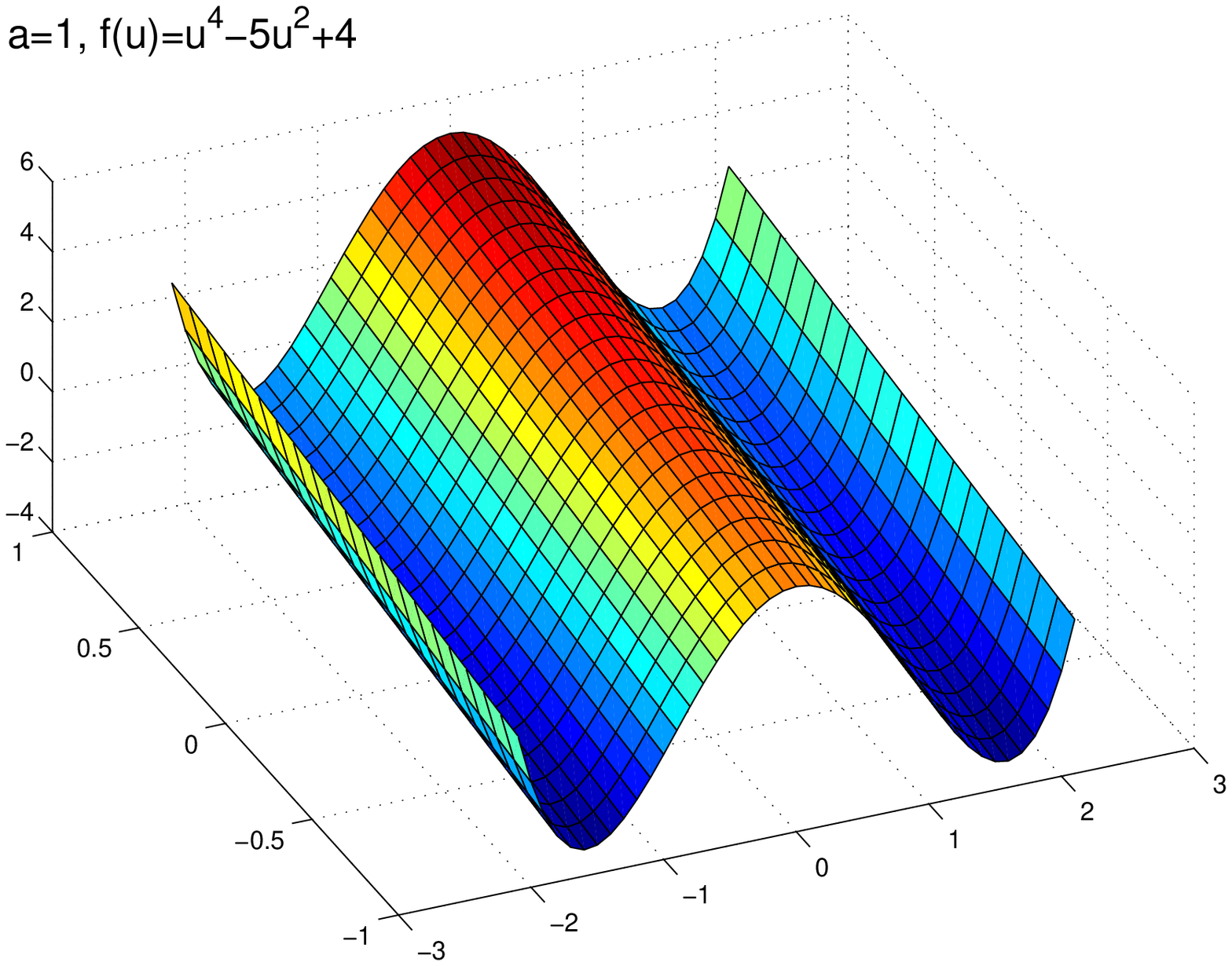}
}
\end{center}
\caption{Cylinders as WPT surfaces}
\end{figure}

\pagebreak

It is interesting to remark that in the case when $f(u)=au+b$ (or $g(v)=av+b$)
(with $a$ and $b$ real constants) and the other function is not polynomial
we still obtain a cylinder, hence a flat surface.
We give in the following two examples:

\begin{figure}[htb]
\begin{center}
\epsfxsize=50mm
\centerline{\leavevmode
\epsffile{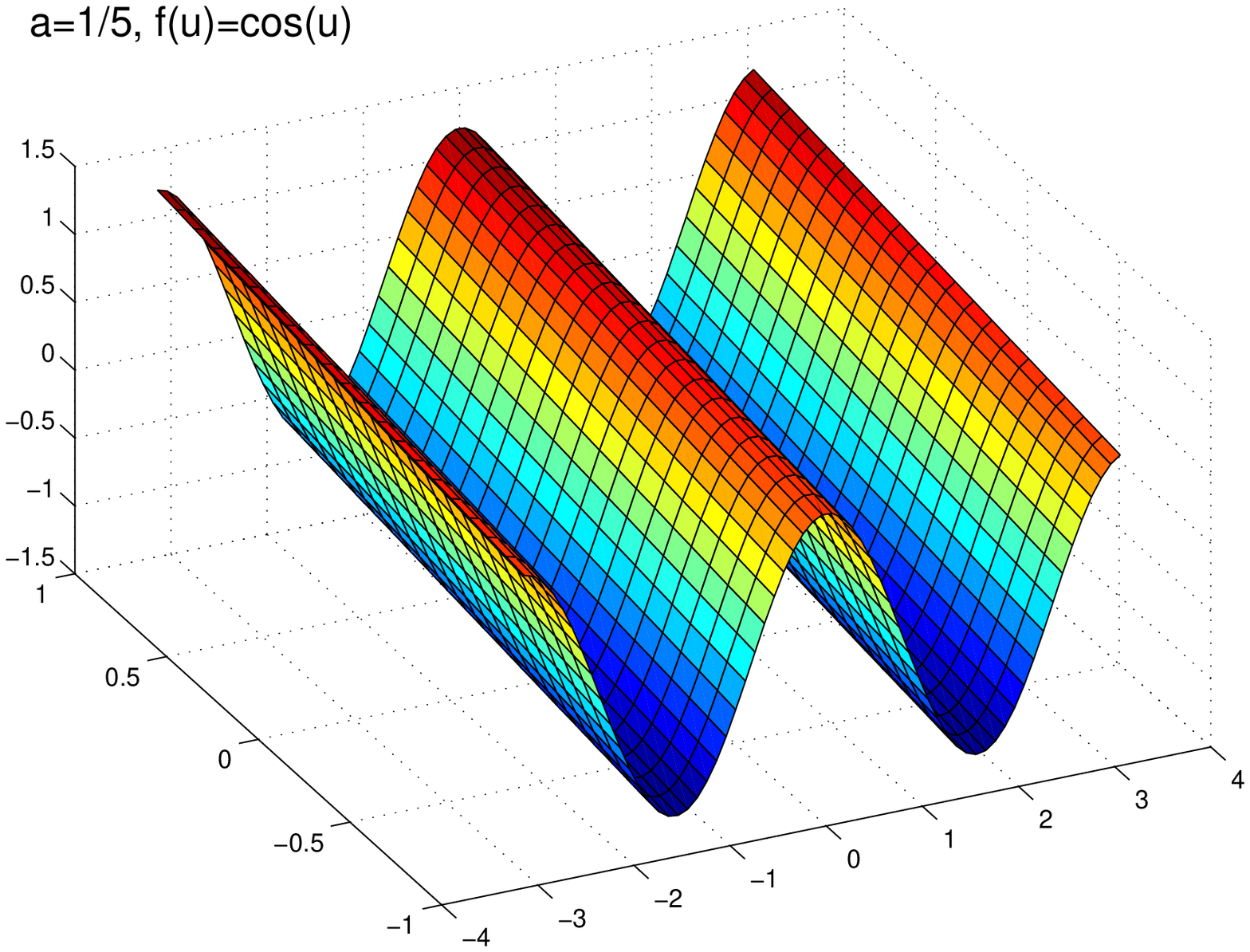}
\qquad
\epsfxsize=50mm
\epsffile{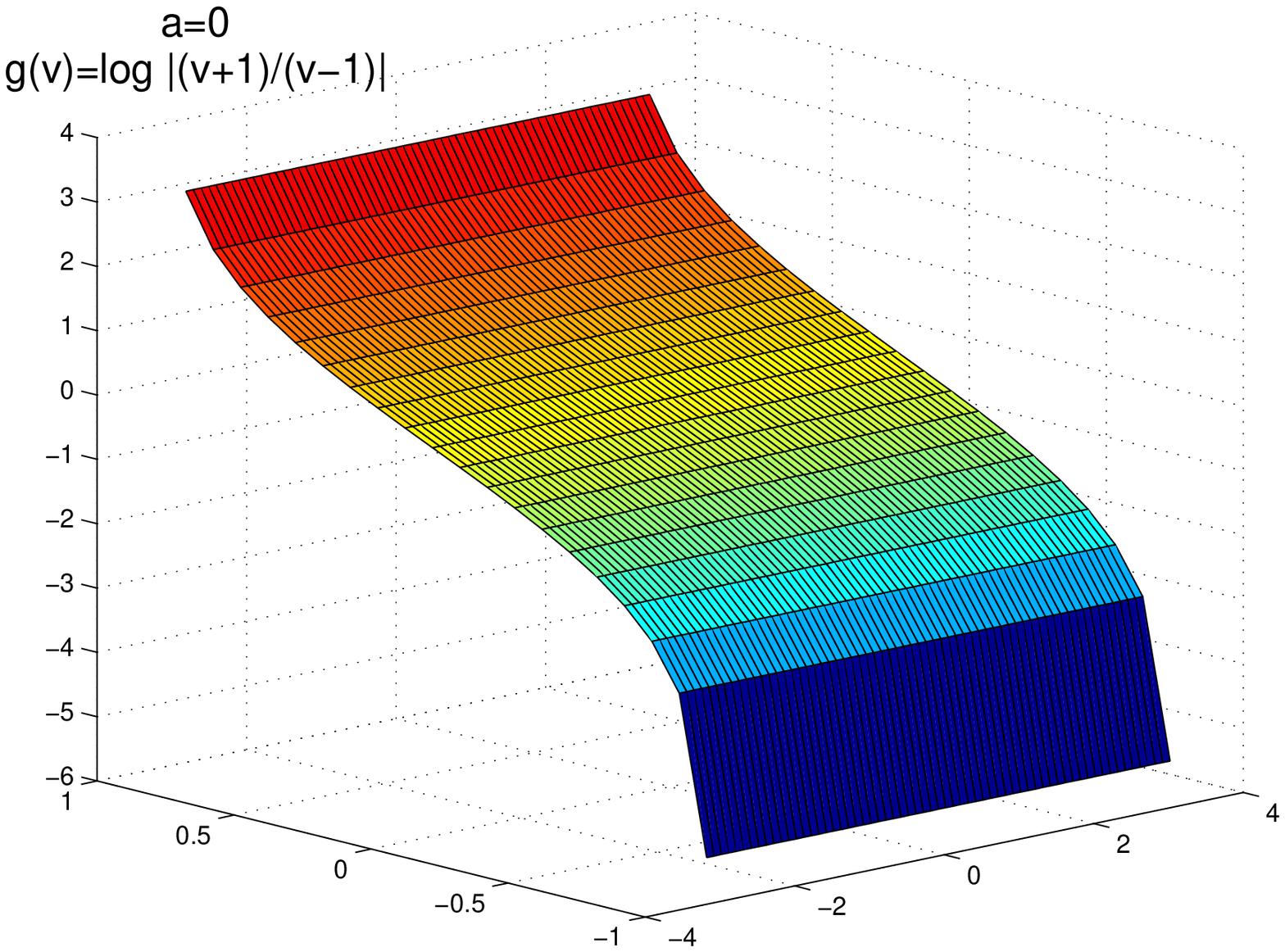}
}
\end{center}
\caption{Other cylinders}
\end{figure}

\section{Linear Weingarten translation surfaces}

A surface $S$ is said to be a linear Weingarten surface -- {\em LW surface} --
if its mean curvature $H$ and Gaussian curvature $K$ satisfy the relation
\begin{equation}
2a H + b K = c
\end{equation}
on $S$ for real numbers $a$, $b$ and $c$ (not all zero). Such a surface is said to be
elliptic, hyperbolic or parabolic depending on whether the discriminant $a^2+bc$
is positive, negative or zero.

\smallskip

Let us analyze the case when $c=0$.

\smallskip

\begin{theorem}
Let $S$ be a linear Weingarten translation surface satisfying
\begin{equation}
\label{eq:lin_c0}
 2a H+bK=0
\end{equation}
$a$ and $b$ real constants, at least one different from zero.
Then
\begin{itemize}
\item [(1)] $S$ is minimal $(H=0)$ {\rm (Scherk's surface)} or
\item [(2)] $S$ is flat $(K=0)$ {\rm (cylinders and planes)}.
\end{itemize}
\end{theorem}
\proof
Consider the parametrization \eqref{eq:param} and the notations
from the previous section. With $H$ and $K$ given by \eqref{eq:H}
and \eqref{eq:K} respectively, we rewrite \eqref{eq:lin_c0} in terms of
$\alpha$ and $\beta$. One gets
$$
a\left( (1+\beta(v)^2)\alpha'(u)+(1+\alpha(u)^2)\beta'(v) \right)+
    \frac{b\alpha'(u)\beta'(v)}{\sqrt{\Delta}}=0
$$
where $\Delta=1+\alpha(u)^2+\beta(v)^2$.

Suppose there exist $u_0$ and $v_0$ such that $\alpha'(u_0)\neq 0$
and $\beta'(v_0)\neq 0$. It follows, due the smoothness of the two functions,
there exist intervals $I\ni u_0$ and $J\ni v_0$ such that $\alpha'\neq0$ on $I$ and
$\beta'\neq0$ on $J$. We obtain
$$
a\left(\frac{1+\alpha(u)^2}{\alpha'(u)}+\frac{1+\beta(v)^2}{\beta'(v)}\right)=
   -\frac b{\sqrt{\Delta}}
$$
for all $u\in I$ and $v\in J$.
The left hand of the above expression is a sum of two terms,
the first depends on $u$ and the second depends on $v$.
Thus, after a derivation with respect to $u$, followed by a
derivation with respect to $v$, the left part of the obtained expression vanishes.
In the same time, the right hand of the expression becomes
$-3b\alpha(u)\beta(v)\alpha'(u)\beta'(v)\Delta^{-\frac 52}$.

Hence $b=0$.

Since $a$ cannot be zero, it follows that $H=0$,
i.e. the surface $S$ is minimal.

The remained cases are $\alpha'=0$ and $\beta'=0$, respectively.
This means that either $f$ or $g$ is a linear function, i.e.
the surface $S$ is a cylinder or a plane.

\endproof

\section{The second Gaussian curvature}

In this section we deal with Riemannian surfaces in Euclidean space $\E^3$
of which the second fundamental form $II$ is positive definite. One may
associate to such surface $S$ geometrical objects measured by means
of its second fundamental form, as {\em second mean curvature} $H_{II}$
and {\em second Gaussian curvature} $K_{II}$, respectively.
Several results concerning $K_{II}$ have been already obtained,
e.g. it is proved that every compact, convex surface in $\E^3$ for which
$K_{II}$ is constant is a sphere \cite{kn:Sch72}.

\medskip

Using the formula of the second Gaussian curvature, a similar one to Brioschi's formula for the Gaussian curvature obtained replacing the
components of the first fundamental form $E$, $F$, $G$ by those of the second fundamental form $e$, $f$, $g$

{\small
$$K_{II}=\frac{1}{(|eg|-f^2)^2}\left(\left|
\begin{array}{l}
-\frac{1}{2}e_{vv}+f_{uv}-\frac{1}{2}g_{uu}\quad \frac{1}{2}e_{u} \quad f_{u}-\frac{1}{2}e_{v}\\[2mm]
f_{v} - \frac{1}{2}g_{u} \hspace{22mm} e \hspace{12mm} f \\[2mm]
\frac{1}{2}g_{v}\hspace{28mm} f \hspace{14mm} g
\end{array}
 \right|
 - \left|
\begin{array}{l}
0 \hspace{6mm} \frac{1}{2}e_{v} \hspace{5mm} \frac{1}{2}g_{u}\\[2mm]
\frac{1}{2}e_{v}\hspace{5mm} e \hspace{7mm} f\\[2mm]
\frac{1}{2}g_{u}\hspace{5mm} f \hspace{7mm} g
\end{array}
\right| \right)$$
}

one gets for the translation surfaces:

$$K_{II}=\frac{num}{4\Delta^{3/2}},$$

where
\begin{equation}
\label{eq:KII}
\begin{array}{c}
num=-2\alpha(u)^2\alpha'(u)^2\beta'(v)-2\alpha'(u)\beta(v)^2\beta'(v)^2+\\
2\alpha(u)^2\alpha'(u)\beta'(v)^2 + 2\alpha'(u)^2\beta(v)^2\beta'(v)+\\
2\alpha'(u)\beta'(v)^2+2\alpha'(u)^2\beta'(v) +\\
\alpha'(u)\beta(v)\beta''(v)+\alpha(u)\alpha''(u)\beta'(v)+\\
+\alpha(u)^2\alpha'(u)\beta(v)\beta''(v)+\alpha(u)\alpha''(u)\beta(v)^2\beta'(v)+\\
 \alpha'(u)\beta(v)^3\beta''(v)+\alpha(u)^3\alpha''(u)\beta'(v).
\end{array}
\end{equation}

\medskip

The general case of the translation surfaces with vanishing second Gaussian curvature in Minkowski 3-space was
recently studied in \cite{kn:Goe07}. An analogous result can be formulated for the Eulidean 3-space.



\medskip

We study the PT surfaces with vanishing second Gaussian curvature, $K_{II}=0$. It turns out that we need to find those
polynomials $\alpha$ and $\beta$ of degree $m$, respectively $n$, so that $num = 0$ in \eqref{eq:KII}.

\medskip
At this point, using the same technique as in {\em Section $2$}  the only convenient case is obtained for $m\geq n=0$, or $n\geq m=0$ namely for
$\beta$ constant or $\alpha$ constant.
We retrieve the {\bf Case 4.} from {\em Preliminaries}, the surfaces given by the parametrization: \eqref{eq:surf_m0} or \eqref{eq:surf_n0}.

\begin{theorem}
The only PT surfaces with vanishing second Gaussian curvature are cylindrical surfaces, which are also WPT surfaces.
\end{theorem}

\medskip

Until now we studied PT surfaces. Inspired from the example given by Blair in \cite{kn:BK92} and mentioned in {\em Preliminaries},
we analyze other types of translation surfaces, involving {\em power functions}, i.e.
\begin{center}
$\alpha = au^p$ and $\beta=bv^q$ with $a,b\in \mathbb{R},\ a,b\neq 0$ and $p,q\in \mathbb{Q}$.
\end{center}

The condition for vanishing second Gaussian curvature becomes:

\begin{equation}
\label{eq:KII_pq}
\begin{array}{c}
a(3p-1)u^pv+a^2b(3q-1)u^{2p+1}v^q + a^3(-p-1)u^{3p}v+\\
b(3q-1)uv^q+ab^2(3p-1)u^pv^{2q+1}+b^3(-q-1)uv^{3q}=0
\end{array}
\end{equation}

\medskip
Again, using the same technique proposed in {\em Section $2$}, the only possible case is for degrees $p=q=\frac{1}{3}$, which implies the additional
condition for coefficients, $a=-b$. We get the surface $S$ given by the parametrization
\begin{equation}
\label{eq:surf_pq1}
r(u,v)=\left( u,v,c(u^{\frac 43}-v^{\frac 43})\right), c\in\R.
\end{equation}
We remark that, up to the multiplication factor $c$, the example given by Blair is the only one
translation surface of this type with vanishing second Gaussian curvature.

\section{Back to the Jacobian condition}

In this section we rewrite the Jacobian condition \eqref{eq:Jac_cond}. Substituting $\alpha$ and $\beta$ with the corresponding power functions,
one gets:

$$
ab^2pq(-q+q^2+pq-pq^2)u^{p+1}v^{2q}+
a^3b^2pq(-2q-pq+2q^2+pq^2)u^{3p+1}v^{2q}+
$$
$$
+a^5b^2pq(-q-2pq+q^2+2pq^2)u^{5p+1}v^{2q}+
ab^4pq(-q+pq-2q^2+2pq^2)u^{p+1}v^{4q}+
$$
\begin{equation}
\label{eq:Jac_cond power}
\begin{array}{c}
+a^3b^4pq(-q-pq-2q^2-4pq^2)u^{3p+1}v^{4q}+
\end{array}
\end{equation}
$$+a^2bpq(p-p^2-pq+p^2q)u^{2p}v^{q+1}+
a^2b^3pq(2p-2p^2+pq-p^2q)u^{2p}v^{3q+1}+
$$
$$
+a^2b^5pq(p-p^2+2pq-2p^2q)u^{2p}v^{5q+1}+
a^4bpq(p+2p^2-pq-2p^2q)u^{4p}v^{q+1}+
$$
$$
+a^4b^3pq(p+2p^2+pq+4p^2q)u^{4p}v^{3q+1}= 0.
$$

It is easy to see that for $p=0$ and any $q\in \mathbb{Q}$ or for $q=0$ and any $p \in \mathbb{Q}$ the identity \eqref{eq:Jac_cond power} is
obviously true. So. if $\alpha=const$ or $\beta=const$, namely $f=au+b$ or $g=cv+d$, with $a,b,c,d\in \mathbb{R}$, then we retrieve the same
parametrization \eqref{eq:surf_m0} or \eqref{eq:surf_n0} like in {\bf Case 4.} from {\em Section $2$}.

\medskip

The other interesting case is $p=q=1$, $\alpha=au$ and $\beta=bv$, retrieving {\bf Case 3.} from {\em Section $2$}, the parametrization of the surface
is given by \eqref{eq:surf_mn1}.

\vspace{5mm}

{\small {\bf Acknowledgement.} The first author was supported by
grant ID\_ 398/2007-2010 , ANCS, Romania. The second author was
partially supported by grant CEEX -- ET n. 5883/2006-2008, ANCS,
Romania. }

\vspace{2mm}

Address (both authors):
\small

University 'Al.I.Cuza' of Ia\c si\\
Faculty of Mathematics\\
Bd. Carol I, no.11\\
700506 Ia\c si\\
Romania\\
e-mail (M.I.Munteanu): marian.ioan.munteanu@gmail.com\\
e-mail (A.I.Nistor): ana.irina.nistor@gmail.com 

\normalsize

\end{document}